\documentclass[journal]{IEEEtran}
\usepackage{epsfig}
\usepackage{graphicx}
\usepackage{psfrag}
\usepackage{ifpdf}
\usepackage{cite}
\usepackage{siunitx}
\usepackage[utf8x]{inputenc}
\usepackage{fancyhdr}
\usepackage{lastpage}
\usepackage{oldstyle}
\usepackage{textcomp}
\usepackage{color}
\usepackage{placeins}
\usepackage{float}
\usepackage{tabularx,colortbl}
\graphicspath{{figures/Appendix/}}
\usepackage{amssymb} % \in
\usepackage{amsmath} % For align command
\usepackage{subfig}
\usepackage{color}
\usepackage{algorithmic}
\usepackage{algorithm}
\usepackage{epstopdf}
\usepackage{morefloats}
\author{Mohammad Tofighi,~\IEEEmembership{Student Member,~IEEE,}
        Onur Yorulmaz
        and~A. Enis Cetin,~\IEEEmembership{Fellow,~IEEE}% <-this % stops a space
\thanks{Authors are with the Department of Electrical and Electronics Engineering, Bilkent University, Ankara, Turkey. M. Tofighi and O. Yorulmaz contributed equally and the names are listed in alphabetical order. emails: tofighi@ee.bilkent.edu.tr, yorulmaz@ee.bilkent.edu.tr, cetin@bilkent.edu.tr.

This work is funded by Turkish Scientific and Technical Research Council (TUBITAK), under project number 113E069.}}

\begin{document}

 \title{Phase and TV Based Convex Sets for Blind Deconvolution of Microscopic Images}

\markboth{Submission to IEEE Journal of Selected Topics in Signal Processing, \today}{}

\maketitle

%\begin{center}
%\normalsize Mohammad Tofighi, Kivanc Kose$^*$, A. Enis Cetin\\
%Dept. of Electrical and Electronic Engineering, Bilkent University,  Ankara, Turkey\\
%$^*$Dermatology Department, Memorial Sloan-Kettering Cancer Center, New York, USA\\
%tofighi@ee.bilkent.edu.tr, $^*$kosek@mskcc.org, cetin@bilkent.edu.tr\\
%\end{center}

\begin{abstract}
%\boldmath
In this article, two closed and convex sets for blind deconvolution problem are proposed. Most blurring functions in microscopy are symmetric with respect to the origin. Therefore, they do not modify the phase of the Fourier transform (FT) of the original image. As a result blurred image and the original image have the same FT phase. Therefore, the set of images with a prescribed FT phase can be used as a constraint set in blind deconvolution problems. Another convex set that can be used during the image reconstruction process is the epigraph set of Total Variation (TV) function. This set does not need a prescribed upper bound on the total variation of the image. The upper bound is automatically adjusted according to the current image of the restoration process. Both of these two closed and convex sets can be used as a part of any blind deconvolution algorithm. Simulation examples are presented.%\footnote{This work is partially presented in IEEE ICIP 2014, Paris, France}
\end{abstract}


\section{Introduction}
\label{sec:Introduction}
A wide range of deconvolution algorithms has been developed to remove blur in microscopic images in recent years \cite{campisi2007blind,Kundur96,Chan98,Sezan91,Sezan199255,Tru84,Xu13,Pengzhao14,Duncan,Sorzano20061205,Acton,Dey,dey:inria-00070726,Pankajakshan:09,Zhang:07}. In this article, two new convex sets are introduced for blind deconvolution algorithms. Both sets can be incorporated to any iterative deconvolution and/or blind deconvolution method.

One of the sets is based on the phase of the Fourier transform (FT) of the observed image.  Most point spread functions blurring microscopic images are symmetric with respect to origin. Therefore, Fourier transform of such functions do not have any phase. As a result, FT phase of the original image and the blurred image have the same phase. The set of images with a prescribed phase is a closed and convex set and projection onto this convex set is easy to perform in Fourier domain.

The second set in the Epigraph Set of Total Variation (ESTV) function. Total variation (TV) value of an image can be limited by an upper-bound to stabilize the restoration process. In fact, such sets were used by many researchers in inverse problems \cite{Com04,Com11,Kos12, PesquetICASSP,ChierchiaSIVP14,dey:inria-00070726}. In this paper, the epigraph of the TV function will be used to automatically estimate an upper-bound on the TV value of a given image. This set is also a closed and convex set. Projection onto ESTV function can be also implemented effectively. ESTV can be incorporated into any iterative blind deconvolution algorithm.

Image reconstruction from Fourier transform phase information was first considered in 1980's \cite{Oppenheim80,Oppenheim81,Oppenheim82,Oppenheim84} and total variation based image denoising was introduced in 1990's \cite{Rud92}. However, FT phase information and ESTV have not been used in blind deconvolution problem to the best of our knowledge.

The paper is organized as follow. In Section II, we review image reconstruction problem from th FT phase and describe the convex set based on phase. In Section III, we describe the epigraph set of TV function.

\section{Convex Set based on the phase of Fourier Transform}
\label{sec:2}
In this section, we introduce our notation and describe how the phase of Fourier transform can be used in deconvolution problems.

Let $x_{o}[n_{1} ,n_{2}]$ be the original image and $h[n_{1} ,n_{2}]$ be the point spread function. The observed image $y$ is obtained by the convolution of $h$ with $x$:
\begin{equation}
\label{eq1}
y[n_{1} ,n_{2}] = h[n_{1} ,n_{2}] \ast x_{o}[n_{1} ,n_{2}],
\end{equation}
where $\ast$ represents the two-dimensional convolution operation. Discrete-time Fourier transform $Y$ of $y$ is, therefore, given by
\begin{equation}
\label{eq2}
Y(w_{1} ,w_{2}) = H(w_{1} ,w_{2})X_{o}(w_{1} ,w_{2}).
\end{equation}
When $h[n_{1} ,n_{2}]$ is symmetric with respect to origin $(h[n_{1} ,n_{2}] = (0, 0))$ phase of $H(w_{1} ,w_{2})$ is zero, i.e., our assumption is $H(w_{1} ,w_{2}) = |H(w_{1} ,w_{2})|$. Point spread functions satisfying this assumption includes uniform Gaussian blurs. Therefore, phase of $Y(w_{1} ,w_{2}) = |Y(w_{1} ,w_{2})|exp(j\measuredangle Y(w_{1} ,w_{2}))$ and $X_{o}(w_{1} ,w_{2}) = |X_{0}(w_{1} ,w_{2})|exp(j\measuredangle X_{o}(w_{1} ,w_{2}))$ are the same:
\begin{equation}
\label{eq3}
\measuredangle Y(w_{1} ,w_{2}) = \measuredangle X_{o}(w_{1} ,w_{2}),
\end{equation}
for all $(w_{1} ,w_{2})$ values. Based on the above observation the following set can be defined:
\begin{equation}
\label{eq4}
C_{\phi} = \{x[n_{1} ,n_{2}]~|~\measuredangle X(w_{1} ,w_{2}) = \measuredangle X_{o}(w_{1} ,w_{2})\},
\end{equation}
which is the set of images whose FT phase is equal to a given prescribed phase $\measuredangle X_{o}(w_{1} ,w_{2})$.

It can easily be shown that this set is closed and convex in $\mathbb{R}^{N_1} \times \mathbb{R}^{N_2}$, for images of size $N_1 \times N_2$.

Projection of an arbitrary image $x$ onto $C_{\phi}$ is implemented in Fourier domain. Let the FT of $x$ be $X(w_{1} ,w_{2}) = |X(w_{1} ,w_{2})| e^{j\phi (w_{1} ,w_{2})}$. The FT $X_p$ of its projection $x_p$ is obtained as follows:
\begin{equation}
\label{eq5}
X_{p}(w_{1} ,w_{2}) = |X(w_{1} ,w_{2})| e^{j\measuredangle X_{o}(w_{1} ,w_{2})},
\end{equation}
where the magnitude of $X_p(w_{1} ,w_{2})$ is the same as the magnitude of $X(w_{1} ,w_{2})$ but its phase is replaced by the prescribed phase function $\measuredangle X_{o}(w_{1} ,w_{2})$. After this step, $x_p [n_{1} ,n_{2}]$ is obtained using the inverse FT. The above operation is implemented using the FFT and implementation details are described in Section \ref{SecIV}.

\begin{figure}[ht]
\centering
\subfloat[]
{\label{fig:Lenanoisy}\includegraphics[width=50mm]{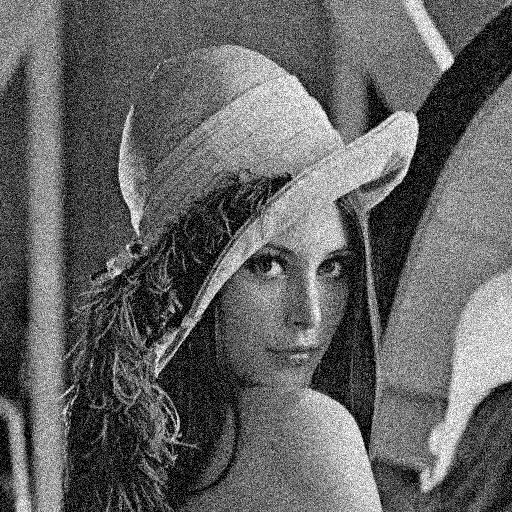}} \quad
\subfloat[]
{\label{fig:Lena_phase}\includegraphics[width=50mm]{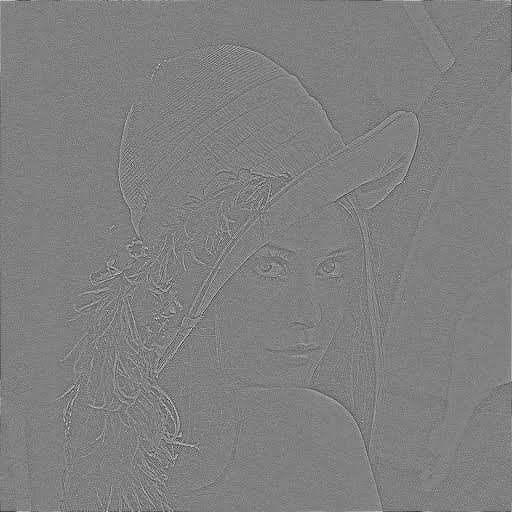}} \quad
\subfloat[]
{\label{fig:Lena_phase_noisy}\includegraphics[width=50mm]{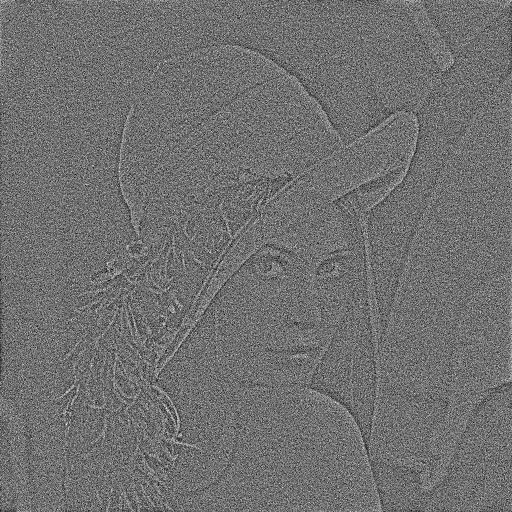}}
\caption{(a) noisy ``Lena" image, (b) Phase only version of the ``Lena" image, and (c) phase only version of the noisy ``Lena" image.}
\label{fig:Lena}
\vspace{-0.5cm}
\end{figure}

Obviously, projection of $y$ onto the set $C_{\phi}$ is the same as itself. Therefore, the iterative blind deconvolution algorithm should not start with the observed image. Image reconstruction from phase (IRP) has been extensively studied by Oppenheim and his coworkers \cite{Oppenheim80,Oppenheim81,Oppenheim82,Oppenheim84}. IRP problem is a robust inverse problem. In Figure \ref{fig:Lena}, phase only version of the well-known Lena image is shown. The phase only image is obtained as follows:
\begin{equation}
\label{eq6}
v = \mathcal{F}^{-1}[C e^{\phi(w_{1} ,w_{2})}]
\end{equation}
where $\mathcal{F}^{-1}{.}$ represents the inverse Fourier transform, $C$ is a constant and $\phi(w_{1} ,w_{2})$ is the phase of Lena image. Edges of the original image are clearly observable in the phase only image. Therefore, the set $C_{\phi}$ contains the crucial edge information of the original image $x_o$.

When the support of $x_o$ is known it is possible to reconstruct the original image from its phase within a scale factor. Oppenheim and coworkers developed Papoulis-Gerchberg type iterative algorithms from a given phase information. In \cite{Oppenheim82} support and phase information are imposed on iterates in space and Fourier domains in a successive manner to reconstruct an image from its phase. 

In blind deconvolution problem the support regions of $x_o$ and $y$ are different from each other. Exact support of the original image is not precisely known; therefore, $C_{\phi}$ is not sufficient by itself to solve the blind deconvolution problem. However, it can be used as a part of any iterative blind deconvolution method.

When there is observation noise, Eq. (\ref{eq1}) becomes:
\begin{equation}
\label{eq7}
\textbf{y}_{o} = \textbf{y} + \boldsymbol{\nu},
\end{equation}
where $\boldsymbol{\nu}$ represents the additive noise. In this case, phase of the observed image is obviously different from the phase of the original image. Luckily, phase information is robust to noise as shown in Fig. \ref{fig:Lena_phase_noisy} which is obtained from a noisy version of Lena image. In spite of noise, edges of Lena are clearly visible in the phase only image. Gaussian noise with variance $\sigma = 30$ is added to Lena image in Fig. \ref{fig:Lenanoisy}.

FTs of some symmetric point spread function may take negative values for some $(w_{1} ,w_{2})$ values. In such $(w_{1} ,w_{2})$ values, phase of the observed image $Y(w_{1} ,w_{2})$ differs from $X(w_{1} ,w_{2})$ by $\pi$. Therefore, phase of $Y(w_{1} ,w_{2})$ should be corrected as in phase unwrapping algorithms. Or some of the $(w_{1} ,w_{2})$ values around $(w_{1},~w_{2}) = (0,~0)$ can be used during the image reconstruction process. It is possible to estimate the main lobe of the FT of the point spread function from the observed image. Phase of FT coefficients within the main lobe are not effected by a shift of $\pi$.

In this article, the set $C_{\phi}$ will be used as a part of the iterative blind deconvolution scheme developed by Dainty \textit{et al} and together with the epigraph set of total variation function which will be introduced in the next section.

\section{Epigraph Set of Total Variation Function}
\label{sec:Convex Minimization}
Bounded total variation is widely used in various image denoising and related applications \cite{PesquetTIP,Chambolle,TofighiICIP,Cen12,Com04,Com11}. The set $C_{\text{TV}}$ of images whose TV values is bounded by a prescribed number $\epsilon$ is defined as follows:
\begin{equation}
\label{eq8}
C_{\text{TV}} = \{{\bf x}: \text{TV}( {\bf x} ) \leq \epsilon \},
\end{equation}
where $\text{TV}$ of an image is defined, in this paper, as follows:
\begin{equation}
\label{eq9}
\text{TV}(\textbf{x}) = \sum_{i,j=1}^{M} |x^{i+1, j} - x^{i, j}| + \sum_{i,j=1}^{M} |x^{i, j+1} - x^{i, j}|.
\end{equation}
This set is closed and convex set in $\mathbb{R}^{N_1 \times N_2}$. Set $C_{\text{TV}}$ can be used in blind deconvolution problems. But the upper bound $\epsilon$ has to be determined somehow a priori. 

In this article we increase the dimension of the space by 1 and consider the problem in $\mathbb{R}^{N_1 \times N_2 + 1}$. We define the epigraph set of the TV function:
\begin{equation}
\label{eq10}
\text{C}_{ESTV} = \{\underline{\mathbf{x}} = [x^T~z]^{T}~|~ \text{TV}(\mathbf{x}) \leq z\},
\end{equation}
where $T$ is the transpose operation and we use bold face letters for $N$ dimensional vectors and underlined bold face letters for $N+1$ dimensional vectors, respectively.

\begin{figure}[ht!]
\begin{center}
\noindent
\includegraphics[width=70mm]{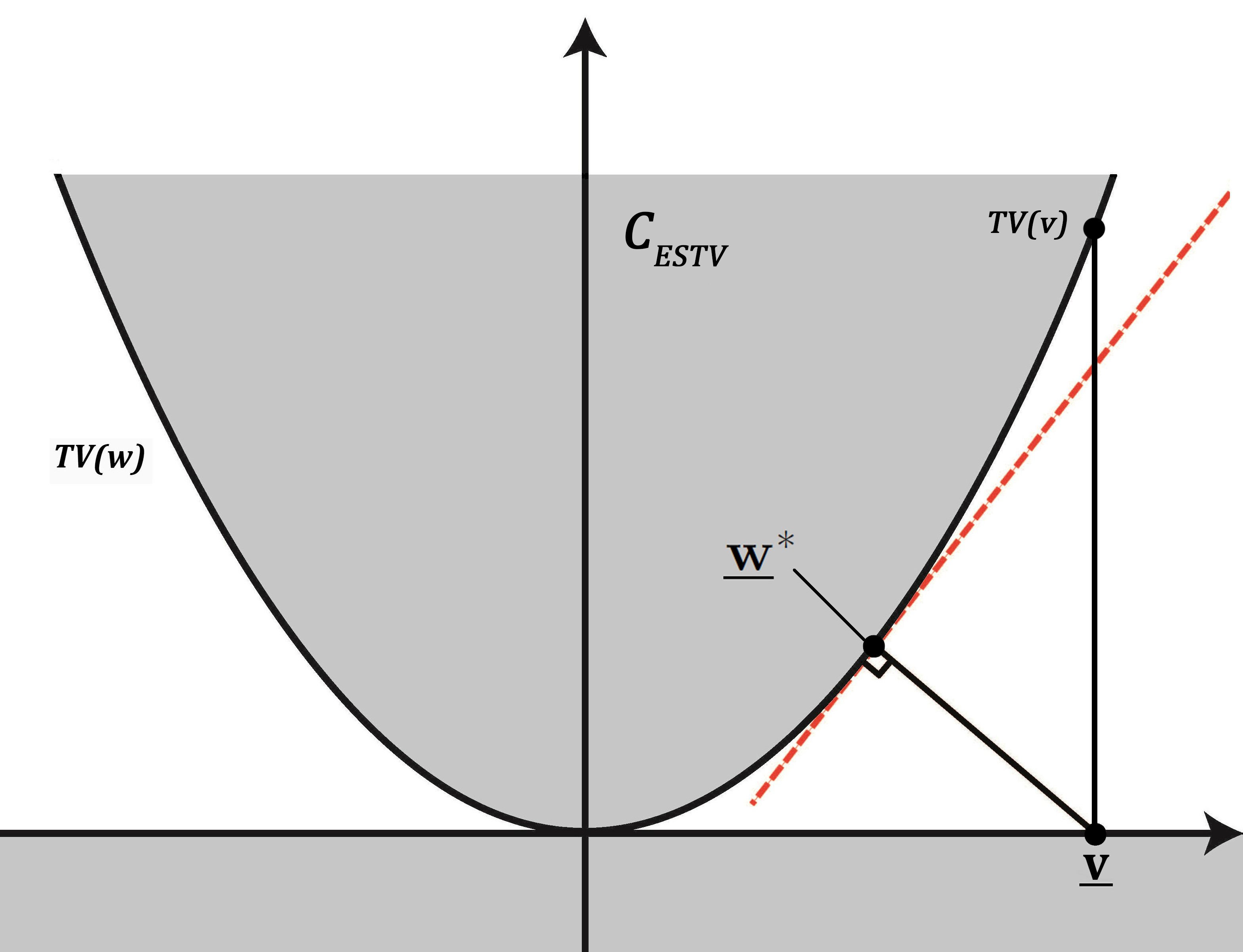}
\caption{Graphical representation of the orthogonal projection onto the set $C_{\text{ESTV}}$ defined in (\ref{eq11}). The observation vector $\underline{\textbf{v}} = [v^T~ 0]^T$ is projected onto the set $C_{\text{ESTV}}$, which is the epigraph set of TV function}
\label{app:convexSol}
\end{center}
\vspace{-0.5cm}
\end{figure}

The concept of the epigraph set is graphically illustrated in Fig. \ref{app:convexSol}. Since $\text{TV}(\mathbf{x})$ is a convex function in $\mathbb{R}^{N_1 \times N_2}$ set the $C_{ESTV}$ is closed and convex in $\mathbb{R}^{N_1 \times N_2 + 1}$. In Eq. (\ref{eq10}) one does not need to specify a prescribed upper bound on TV of an image. An orthogonal projection onto the set $C_{ESTV}$ reduces the total variation value of the image as graphically illustrated in Fig. \ref{app:convexSol} because of the convex nature of the TV function. Let $\textbf{v}$ be an $N = N_{1}\times N_{2}$ dimensional image to be projected onto the set $C_{ESTV}$. In orthogonal projection operation, we select the nearest vector $\underline{\mathbf{x}}^{\star}$ on the set $C_{ESTV}$ to $\underline{\mathbf{w}}$. The projection vector $\mathbf{x}^{\star}$ of an image $\textbf{v}$ is defined as:
\begin{equation}
\label{eq11}
\underline{\mathbf{w}}^{\star} = \text{arg} \underset{\underline{\mathbf{w}}\in C_{\text{ESTV}}}{\text{min}} \|\underline{\mathbf{v}} - \underline{\mathbf{w}}\|^{2},
\end{equation}
where $\underline{\mathbf{v}}$ = [$\textbf{v}^{T}$\ 0]. The projection operation described in (\ref{eq11}) is equivalent to:
\begin{equation}
\underline{\mathbf{w}}^{\star} = \begin{bmatrix} \textbf{w}_p \\ \text{TV}(\textbf{w}_p) \end{bmatrix} = \text{arg} \underset{\underline{\mathbf{w}}\in \text{C}_{\mathrm{f}}}{\text{min}} \|\begin{bmatrix} \textbf{v} \\ 0 \end{bmatrix} - \begin{bmatrix} \textbf{w} \\ \text{TV}(\textbf{w}) \end{bmatrix}\|,
\label{n1}
\end{equation}
where $\underline{\mathbf{w}}^{\star}=[\textbf{w}_p^T, \text{TV}(\textbf{w}_p)]$ is the projection of $(\textbf{v},0)$ onto the epigraph set.  The projection $\underline{\mathbf{w}}^{\star}$ must be on the boundary of the epigraph set. Therefore, the projection must be on the form $[\textbf{w}_p^T, \text{TV}(\textbf{w}_p)]$. Equation (\ref{n1}) becomes:
\begin{equation}
\underline{\mathbf{w}}^{\star} = \begin{bmatrix} \textbf{w}_p \\ \text{TV}(\textbf{w}_p) \end{bmatrix} = \text{arg} \underset{\underline{\mathbf{w}}\in \text{C}_{\mathrm{f}}}{\text{min}}\|\textbf{v}-\textbf{w}\|_2^2 + \text{TV}(\textbf{w})^2.
\label{n2}
\end{equation}
It is also possible to use $\lambda\text{TV}(.)$ as a the convex cost function and Eq. \ref{n2} becomes:

\begin{equation}
\underline{\mathbf{w}}^{\star} = \begin{bmatrix} \textbf{w}_p \\ \text{TV}(\textbf{w}_p) \end{bmatrix} = \text{arg} \underset{\underline{\mathbf{w}}\in \text{C}_{\mathrm{f}}}{\text{min}}\|\textbf{v}-\textbf{w}\|_2^2 + \lambda^2\text{TV}(\textbf{w})^2.
\label{n3}
\end{equation}
The solution of (\ref{eq11}) can be obtained using the method that we discussed in \cite{TofighiICIP,GlobalSIP2013}. The solution is obtained in an iterative manner and the key step in each iteration is an orthogonal projection onto a supporting hyperplane of the set $C_{ESTV}$.

In current TV based denoising methods \cite{Chambolle, Com11} the following cost function is used:
\begin{equation}
\label{app:eq:cost}
{{\text{min}} \|\mathbf{v}} - {\mathbf{w}}\|^2_2 + \lambda \text{TV}(\textbf{w}).
\end{equation}
However, we were not able to prove that \ref{app:eq:cost} corresponds to a non-expensive map or not. On the other hand, minimization problem in Eq. (\ref{n2}) and (\ref{n3}) are the results of projection onto convex sets, as a result they correspond to non-expensive maps \cite{You82,Cetin1990,Cetin88,KivancSPM,PesquetTIP,Com04,Tru85,You82,Stark:81,Sezan199255,Combettes93}. Therefore, they can be incorporated into any iterative deblurring algorithm without effecting the convergence of the algorithm.

\section{How to incorporate $C_{ESTV}$ and $C_{\phi}$ into a deblurring method}
\label{SecIV}
One of the earliest blind deconvolution methods is the iterative space-Fourier domain method developed by Ayers and Dainty \cite{Ayers:88}. In this approach, iterations start with a $x_{o}[n] = x_{o}[n_{1}, n_{2}]$, where we introduce a new notation to specify equations $[n] = [n_{1}, n_{2}]$. For example, we rewrite Eq. (\ref{eq1}) as follows:
\begin{equation}
\label{eq12}
y[n] = h[n]\ast x_{o}[n],
\end{equation}
The method successively updates $h[n]$ and $x[n]$ in a Wiener filter-like equation. Here is the $i^{th}$ step of the algorithm:
\begin{enumerate}
\item Compute $\hat{X}_{i}(w) = \mathcal{F}\{x_{i}[n]\}$, where $\mathcal{F}$ represents the FT operation and $w = (w_{1}, w_{2})$, with some abuse of notation.
\item Estimate the point-spread filter response using the following equation
\begin{equation}
\label{eq13}
\tilde{H}_{i}(w) = \frac{Y(w)\hat{X}_{i}^{*}(w)}{|\hat{X}_{i}(w)|^{2} + \alpha /|\hat{H}_{i}(w)|^{2}},
\end{equation}
where $\alpha$ is a small real number.
\item Compute $\tilde{h}_{i}[n] = \mathcal{F}^{-1}\{\tilde{H}_{i}(w)\}$
\item Impose the positivity constraint and finite support constraints on $\tilde{h}_{i}[n]$. Let the output of this step be $\hat{h}_{i}[n]$.
\item Compute $\hat{H}_{i}(w) = \mathcal{F}\{\hat{h}_{i}[n]\}$
\item Update the image
\begin{equation}
\label{eq14}
\tilde{X}_{i}(w) = \frac{Y(w)\hat{H}_{i}^{*}(w)}{|\hat{H}_{i}(w)|^{2} + \alpha /|\hat{X}_{i}(w)|^{2}},
\end{equation}
\item Compute $\hat{x}_{i}[n] = \mathcal{F}^{-1}\{\hat{X}_{i}(w)\}$
\item Impose spatial domain positivity and finite support constraint on $\hat{x}_{i}[n]$ to produce the next iterate $\hat{x}_{i+1}[n]$.
\end{enumerate}

\begin{figure}[ht!]
\begin{center}
\noindent
\includegraphics[width=85mm]{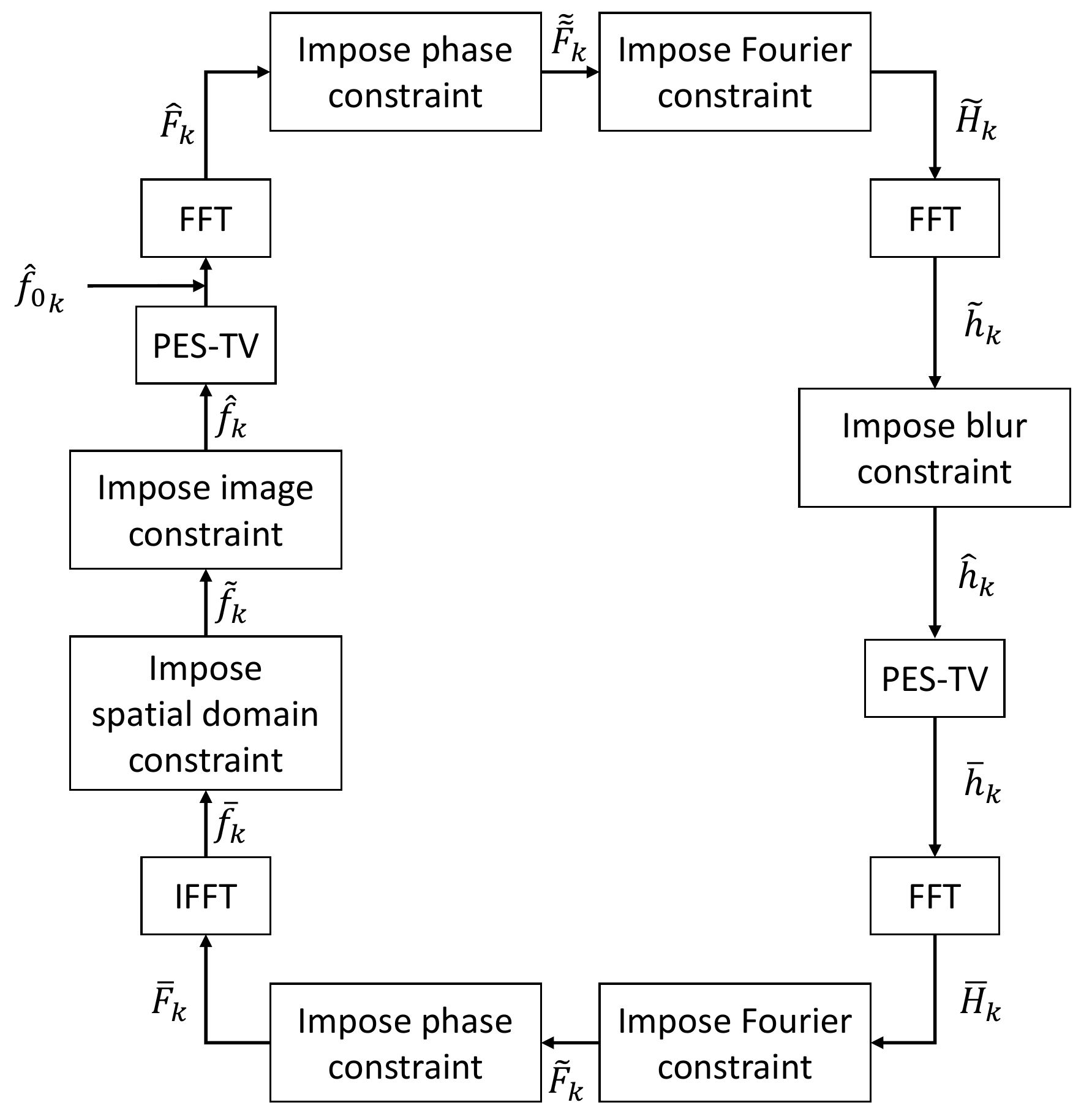}
\caption{Flow chart of the proposed algorithm. PES-TV stands for Projection onto the Epigraph Set of TV function.}
\label{app:Flow}
\end{center}
\vspace{-0.5cm}
\end{figure}

Iterations are stopped when there is no significant change between successive iterates. We can easily modify this algorithm using the convex sets defined in Section \ref{sec:2} and \ref{sec:Convex Minimization}. We introduce step 6-a as follows:

6-a) Impose the phase information
\begin{equation}
\label{eq15}
\bar{X}_{i}(w) = |\tilde{X}_{i}(w)| e^{j\measuredangle Y(w)},
\end{equation}
where $\measuredangle Y(w)$ is the phase of $Y(w)$. This step is the projection onto the set $C_{\phi}$. As a result step 7 becomes $\tilde{x}_{i}[n] = \mathcal{F}^{-1}\{\bar{X}_{i}(w)\}$. We also introduce a new step to Ayers and Dainty's algorithm as follows: Project $\tilde{x}_{i}[n]$ onto the set $C_{\text{ESTV}}$ to obtain $\hat{x}_{i+1}[n]$. The flowchart of the proposed algorithm is shown in Fig. \ref{app:Flow}.

Since the filter is a zero-phase filter in microscompic image analysis $h[n_1, n_2] = h[-n_1, -n_2]=h[-n_1, n_2]=h[n_1, -n_2]$ this condition is also imposed on the current iterate in Step 4. 

Global convergence of Ayers-Dainty method has not been proved. In fact, we experimentally observed that it may diverge in some FL microscopy images. Projections onto convex sets are non-expansive maps \cite{Combettes93,Slavakis,Yamada2001473}, therefore, they do not cause any divergence problems in an iterative image debluring algorithm.

\begin{table*}[htbp]
  \centering
  \caption{Deconvolution results for florescence microscopic images blurred by Gaussian filter with disc size $d$ and $\sigma = 1$. PSNR (dB) values are presented for the proposed algorithm and Ayers and Dainty's algorithm.}
    \label{Tab1}
    \begin{tabular}{|c|c|c|c|c|c|c|c|c|c|c|c|c|c|c|c|c|}
	\hline
    Method & Filter radius & Im-1  & Im-2  & Im-3  & Im-4  &  Im-5 & Im-6  & Im-7  & Im-8  & Im-9  & Im-10 & Im-11 & Im-12 & Im-13 & Im-14 \\\hline\hline
    Ayers & d = 5 & 5.39  & \textbf{17.98} & 6.13  & 9.20  & 8.22  & 5.52  & 12.31 & 7.52  & 6.79  & 5.67  & 6.53  & 13.34 & 14.49 & 8.13 \\\hline
        Modified & d = 5 & \textbf{9.16} & 11.00 & \textbf{10.36} & \textbf{11.80} & \textbf{9.85} & \textbf{10.37} & \textbf{16.49} & \textbf{8.90} & \textbf{10.86} & \textbf{8.47} & \textbf{10.34} & \textbf{16.72} & \textbf{17.48} & \textbf{10.68} \\\hline
        Ayers & d = 10 & 4.57  & 5.25  & 5.12  & 8.70  & 7.61  & 10.93 & 13.51 & 6.52  & 7.70  & 5.96  & 6.44  & 10.17 & 16.80 & 7.53 \\\hline
        Modified & d = 10 & \textbf{11.14} & \textbf{9.41} & \textbf{9.81} & \textbf{11.61} & \textbf{9.96} & \textbf{12.80} & \textbf{16.04} & \textbf{9.45} & \textbf{11.45} & \textbf{16.43} & \textbf{11.71} & \textbf{15.42} & \textbf{18.14} & \textbf{11.37} \\\hline
        Ayers & d = 15 & 4.94  & 5.30  & 5.90  & 7.81  & 6.33  & 7.76  & 11.49 & 8.47  & 5.56  & 4.01  & 6.38  & 10.23 & \textbf{18.67} & 7.22 \\\hline
        Modified & d = 15 & \textbf{9.00} & \textbf{10.60} & \textbf{11.51} & \textbf{12.04} & \textbf{10.31} & \textbf{10.49} & \textbf{14.16} & \textbf{10.72} & \textbf{10.99} & \textbf{7.83} & \textbf{9.49} & \textbf{15.65} & 17.69 & \textbf{10.71} \\\hline
    \end{tabular}%
\end{table*}%

\begin{table*}[htbp]
  \centering
  \caption{Deconvolution results for florescence microscopic images blurred by Gaussian filter with disc size $d$ and $\sigma = 2$. PSNR (dB) values are presented for the proposed algorithm and Ayers and Dainty's algorithm.}
    \label{Tab2}
    \begin{tabular}{|c|c|c|c|c|c|c|c|c|c|c|c|c|c|c|c|c|}
	\hline
	Method & Filter radius & Im-1  & Im-2  & Im-3  & Im-4  &  Im-5 & Im-6  & Im-7  & Im-8  & Im-9  & Im-10 & Im-11 & Im-12 & Im-13 & Im-14 \\\hline\hline
    Ayers & d = 5 & 6.37  & 6.92  & 7.03  & 8.10  & 9.08  & 8.69  & 15.64 & 8.44  & \textbf{10.17} & \textbf{14.19} & \textbf{10.98} & 16.28 & \textbf{22.95} & 11.41 \\\hline
    Modified & d = 5 & \textbf{17.24} & \textbf{12.61} & \textbf{11.79} & \textbf{12.08} & \textbf{12.82} & \textbf{10.36} & \textbf{18.81} & \textbf{12.36} & 9.38  & 13.20 & 8.40  & \textbf{16.66} & 17.07 & \textbf{13.27} \\\hline
    Ayers & d = 10 & \textbf{16.36} & 7.04  & 8.04  & 9.63  & \textbf{24.20} & 10.58 & 19.26 & 9.58  & \textbf{16.78} & 8.07  & 6.57  & 18.63 & \textbf{22.35} & 15.21 \\\hline
    Modified & d = 10 & 15.30 & \textbf{13.01} & \textbf{12.58} & \textbf{11.86} & 16.86 & \textbf{16.33} & \textbf{20.44} & \textbf{11.61} & 10.47 & \textbf{14.18} & \textbf{11.74} & \textbf{22.17} & 21.66 & \textbf{18.39} \\\hline
    Ayers & d = 15 & 10.52 & \textbf{13.85} & 12.71 & \textbf{14.85} & 13.49 & 15.64 & 18.38 & 9.01  & 7.20  & 7.10  & 6.35  & \textbf{22.07} & 19.60 & 13.48 \\\hline
    Modified & d = 15 & \textbf{20.96} & 11.70 & \textbf{16.14} & 12.99 & \textbf{20.76} & \textbf{19.80} & \textbf{21.23} & \textbf{15.08} & \textbf{11.42} & \textbf{13.84} & \textbf{12.32} & 21.91 & \textbf{23.05} & \textbf{18.12} \\\hline
    \end{tabular}%
\end{table*}%

\begin{table*}[htbp]
  \centering
  \caption{Deconvolution results for florescence microscopic images blurred by Gaussian filter with disc size $d$ and $\sigma = 3$. PSNR (dB) values are presented for the proposed algorithm and Ayers and Dainty's algorithm.}
  \label{Tab3}
    \begin{tabular}{|c|c|c|c|c|c|c|c|c|c|c|c|c|c|c|c|c|}
	\hline
    Method & Filter radius & Im-1  & Im-2  & Im-3  & Im-4  &  Im-5 & Im-6  & Im-7  & Im-8  & Im-9  & Im-10 & Im-11 & Im-12 & Im-13 & Im-14 \\\hline\hline
%    \midrule
    Ayers & d = 5 & 7.40  & 6.33  & 6.39  & 7.43  & 6.51  & 9.39  & 15.92 & 7.20  & 6.35  & 6.71  & 6.62  & 17.47 & 22.19 & 8.88 \\\hline
    Modified & d = 5 & \textbf{8.08} & \textbf{23.06} & \textbf{17.18} & \textbf{22.16} & \textbf{11.66} & \textbf{23.18} & \textbf{21.75} & \textbf{23.83} & \textbf{20.26} & \textbf{23.20} & \textbf{30.91} & \textbf{17.77} & \textbf{23.59} & \textbf{8.91} \\\hline
    Ayers & d = 10 & \textbf{8.80} & 20.03 & 14.95 & 16.08 & 22.06 & 22.57 & 21.18 & 22.97 & 17.39 & 22.03 & 32.26 & 21.55 & 24.74 & 20.89 \\\hline
    Modified & d = 10 & 8.02  & \textbf{25.66} & \textbf{24.74} & \textbf{26.44} & \textbf{24.15} & \textbf{29.61} & \textbf{23.99} & \textbf{24.44} & \textbf{20.67} & \textbf{26.04} & \textbf{39.92} & \textbf{24.03} & \textbf{27.05} & \textbf{27.30} \\\hline
    Ayers & d = 15 & 18.29 & \textbf{28.88} & 14.36 & 18.77 & 27.50 & 27.54 & 24.15 & \textbf{24.91} & 21.64 & 26.09 & 34.00 & 22.60 & 23.69 & 26.76 \\\hline
    Modified & d = 15 & \textbf{23.86} & 28.62 & \textbf{32.23} & \textbf{28.55} & \textbf{36.93} & \textbf{27.80} & \textbf{24.31} & 24.33 & \textbf{21.34} & \textbf{29.63} & \textbf{40.84} & \textbf{23.16} & \textbf{27.44} & \textbf{35.31} \\\hline
    \end{tabular}%
\end{table*}%

\section{Experimental Results}
The proposed algorithm is evaluated using different florescence (FL) microscopy images obtained at Bilkent University. Gaussian and uniform filters are used to blur the images. These images are blurred by Gaussian filter with disc sizes $d = 5, 10,~\text{and}~15$ and $\sigma = 1, 2,~\text{and}~3$. The blind deconvolution results are presented for $\sigma = 1, 2,~\text{and}~3$ in Tables \ref{Tab1}, \ref{Tab2}, and \ref{Tab3}, respectively.

In Tables \ref{Tab1}, \ref{Tab2}, and \ref{Tab3} bold font is used for the highest PSNR. Clearly, the modified deblurring method using $C_{\phi}$ and $C_{\text{ESTV}}$ produces better PSNR values in almost all cases. 

\begin{figure}[ht]
\centering
\subfloat[]
{\label{fig:F1}\includegraphics[width=40mm]{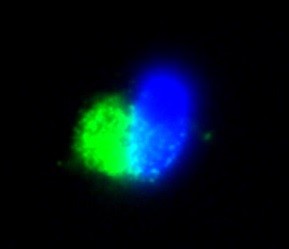}} \quad
\subfloat[]
{\label{fig:F2}\includegraphics[width=40mm]{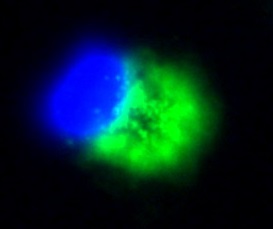}} \\
\subfloat[]
{\label{fig:F3}\includegraphics[width=40mm]{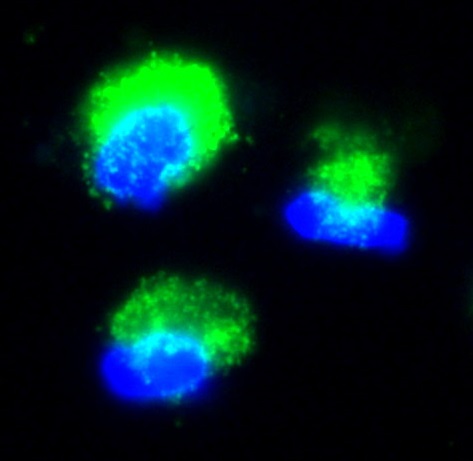}} \quad
\subfloat[]
{\label{fig:F4}\includegraphics[width=40mm]{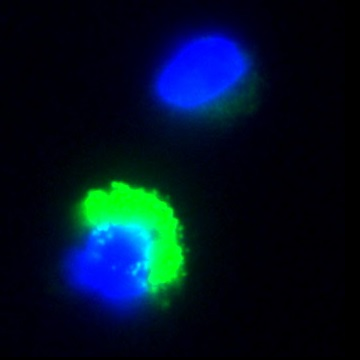}}
\caption{Sample  florescence microscopic images used in experiments (a) Im-5, (b) Im-7 (c) Im-10, and (d) Im-11.}
\label{fig:sample}
\end{figure}

\begin{figure*}[ht]
\centering
\subfloat[]
{\label{fig:origIm11}\includegraphics[width=58mm]{a11.png}} \quad
\subfloat[]
{\label{fig:BlurredIm11}\includegraphics[width=58mm]{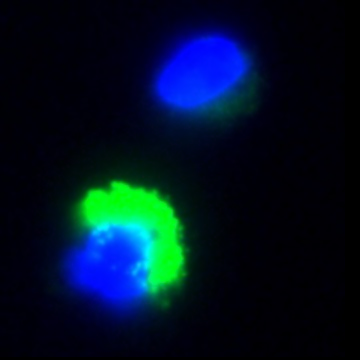}} \\
\subfloat[]
{\label{fig:ModAyersDebIm11}\includegraphics[width=58mm]{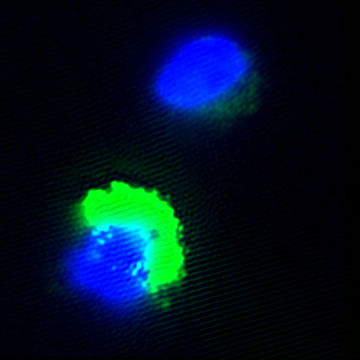}} \quad
\subfloat[]
{\label{fig:AyersDebIm11}\includegraphics[width=58mm]{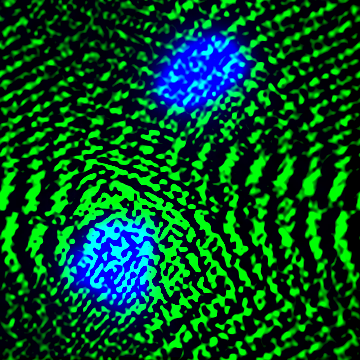}}
\caption{The deconvolution results for Im-11 (a) Original Im-11, (b) blurred Im-11 ($d = 5$, $\sigma = 3$) (c) Deblurred by the proposed algorithm (PSNR = 30.91 dB), and (d) Deblurred by Ayers and Dainty's algorithm (PSNR = 6.62 dB).}
\label{fig:a11}
\end{figure*}

Four sample of images used in this set of experiments are shown in Fig. \ref{fig:sample}. As an example, the Im-11 shown in Fig. \ref{fig:origIm11} in is blurred using Gaussian filter with $d = 5$ and $\sigma = 3$. The blurred image is shown in Fig. \ref{fig:BlurredIm11}. The deblurred image obtained using the proposed algorithm and Ayers and Dainty's algorithm are shown in Fig. \ref{fig:ModAyersDebIm11} and \ref{fig:AyersDebIm11}, respectively.

\begin{figure*}[ht]
\centering
\subfloat[]
{\label{fig:origflor}\includegraphics[width=70mm]{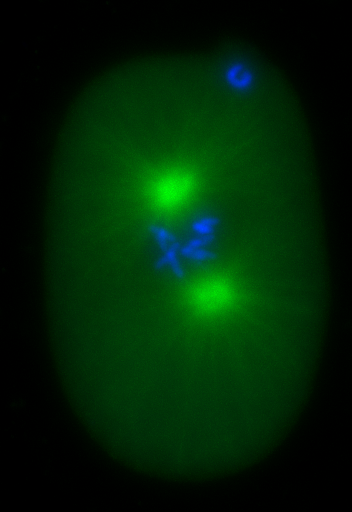}} \quad
\subfloat[]
{\label{fig:Deblurredflor}\includegraphics[width=70mm]{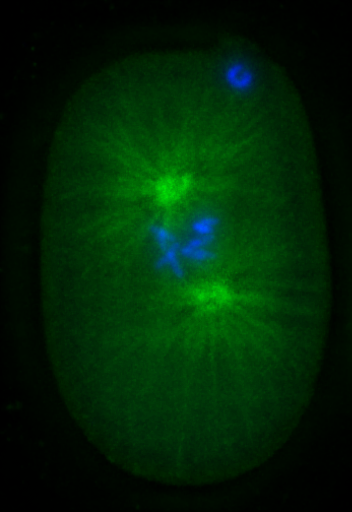}} \\
\subfloat[]
{\label{fig:AyersDebflor}\includegraphics[width=70mm]{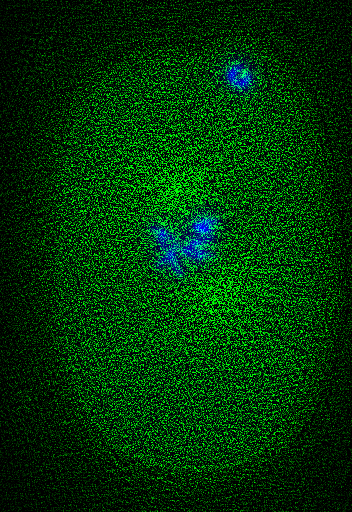}} \quad
\subfloat[]
{\label{fig:ChineDebflor}\includegraphics[width=70mm]{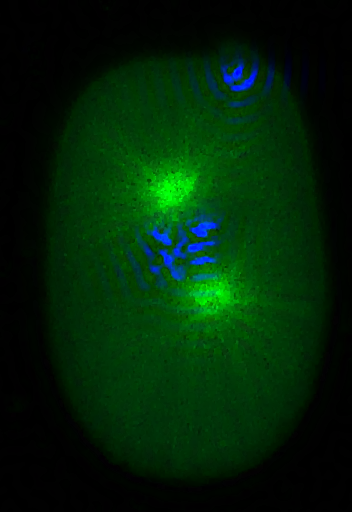}}
\caption{The deconvolution results for FL image downloaded from [http://bigwww.epfl.ch/algorithms/mltldeconvolution/] (a) blurred image, (b) deblurred by the blind deconvolution using phase information, (c) deblurred by Ayers and Dainty’s algorithm (PSNR = 34.00 dB), and (d) Deblurred by Xu et al’s algorithm \cite{Xu13}.}
\label{fig:flor}
\end{figure*}

In another set of experiments, we used the FL image shown in Fig. \ref{fig:origflor} which is blurred by an unknown filter or captured with a focus blur \cite{Vonesch09}. This image is deblurred using the blind deconvolution by phase information and its output is compared with Ayers and Dainty's and Xu \textit{et al}'s algorithm \cite{Xu13}. The deblurred image using the blind deconvolution by phase information and $C_\text{ESTV}$ , Ayers and Dainty’s algorithm, and the Xu \textit{et al}’s algorithm are shown in Fig. \ref{fig:Deblurredflor}, \ref{fig:AyersDebflor}, and \ref{fig:ChineDebflor}, respectively.

Ayers and Dainty's method does not converge as shown in Fig. \ref{fig:AyersDebflor}. Xu \textit{et al}'s algorithm also diverges when we select ``default" option. It does not diverge when we select ``small" kernel option but the result is far from perfect as shown in Fig. \ref{fig:flor}. Sets $C_{\phi}$ and $C_{\text{ESTV}}$ can be also incorporated into
Xu \textit{et al}'s method for symmetric kernels but we do not have an access to the source code. We get the best results when we use $C_{\phi}$ and $C_{\text{ESTV}}$ in a successive manner as shown in Fig. \ref{fig:ModAyersDebIm11} and \ref{fig:Deblurredflor}.

Iterations are stopped after 300 rounds in all cases. In the following web-page you may find the MATLAB code of projections onto $C_{\phi}$ and $C_{\text{ESTV}}$ and the example FL images which four of them are shown in Fig. \ref{fig:sample}. Web-page: http://signal.ee.bilkent.edu.tr/BlindDeconvolution.html.

\section{Conclusion}
FT phase and the epigraph of the TV function are closed and convex sets. They can be used as a part of iterative microscopic image deblurring algorithms. Both sets not only provide additional information about the desired solution but they also stabilize the deconvolution algorithms. We experimentally observed that they significantly improved the deblurring results of Ayers and Dainty's method.
%\clearpage
\bibliographystyle{IEEEtran}
\bibliography{PhdReferences}

\begin{thebibliography}{10}
\providecommand{\url}[1]{#1}
\csname url@samestyle\endcsname
\providecommand{\newblock}{\relax}
\providecommand{\bibinfo}[2]{#2}
\providecommand{\BIBentrySTDinterwordspacing}{\spaceskip=0pt\relax}
\providecommand{\BIBentryALTinterwordstretchfactor}{4}
\providecommand{\BIBentryALTinterwordspacing}{\spaceskip=\fontdimen2\font plus
\BIBentryALTinterwordstretchfactor\fontdimen3\font minus
  \fontdimen4\font\relax}
\providecommand{\BIBforeignlanguage}[2]{{%
\expandafter\ifx\csname l@#1\endcsname\relax
\typeout{** WARNING: IEEEtran.bst: No hyphenation pattern has been}%
\typeout{** loaded for the language `#1'. Using the pattern for}%
\typeout{** the default language instead.}%
\else
\language=\csname l@#1\endcsname
\fi
#2}}
\providecommand{\BIBdecl}{\relax}
\BIBdecl

\bibitem{campisi2007blind}
P.~Campisi and K.~Egiazarian, \emph{Blind image deconvolution: theory and
  applications}.\hskip 1em plus 0.5em minus 0.4em\relax CRC press, 2007.

\bibitem{Kundur96}
D.~Kundur and D.~Hatzinakos, ``Blind image deconvolution,'' \emph{Signal
  Processing Magazine, IEEE}, vol.~13, no.~3, pp. 43--64, May 1996.

\bibitem{Chan98}
T.~Chan and C.-K. Wong, ``Total variation blind deconvolution,'' \emph{Image
  Processing, IEEE Transactions on}, vol.~7, no.~3, pp. 370--375, Mar 1998.

\bibitem{Sezan91}
M.~Sezan and H.~Trussell, ``Prototype image constraints for set-theoretic image
  restoration,'' \emph{Signal Processing, IEEE Transactions on}, vol.~39,
  no.~10, pp. 2275--2285, Oct 1991.

\bibitem{Sezan199255}
M.~Sezan, ``An overview of convex projections theory and its application to
  image recovery problems,'' \emph{Ultramicroscopy}, vol.~40, no.~1, pp. 55 --
  67, 1992.

\bibitem{Tru84}
H.~Trussell and M.~Civanlar, ``The feasible solution in signal restoration,''
  \emph{Acoustics, Speech and Signal Processing, IEEE Transactions on},
  vol.~32, no.~2, pp. 201--212, Apr 1984.

\bibitem{Xu13}
L.~Xu, S.~Zheng, and J.~Jia, ``Unnatural l0 sparse representation for natural
  image deblurring,'' in \emph{Computer Vision and Pattern Recognition (CVPR),
  2013 IEEE Conference on}, June 2013, pp. 1107--1114.

\bibitem{Pengzhao14}
P.~Ye, H.~Feng, Q.~Li, Z.~Xu, and Y.~Chen, ``Blind deconvolution using an
  improved l0 sparse representation,'' pp. 928\,419--928\,419--6, 2014.

\bibitem{Duncan}
J.~Boulanger, C.~Kervrann, and P.~Bouthemy,
  ``\BIBforeignlanguage{English}{Adaptive spatio-temporal restoration for 4d
  fluorescence microscopic imaging},'' in
  \emph{\BIBforeignlanguage{English}{Medical Image Computing and
  Computer-Assisted Intervention – MICCAI 2005}}, ser. Lecture Notes in
  Computer Science, J.~Duncan and G.~Gerig, Eds.\hskip 1em plus 0.5em minus
  0.4em\relax Springer Berlin Heidelberg, 2005, vol. 3749, pp. 893--901.

\bibitem{Sorzano20061205}
C.~Sorzano, E.~Ortiz, M.~López, and J.~Rodrigo, ``Improved bayesian image
  denoising based on wavelets with applications to electron microscopy,''
  \emph{Pattern Recognition}, vol.~39, no.~6, pp. 1205 -- 1213, 2006.

\bibitem{Acton}
S.~Acton, ``Deconvolutional speckle reducing anisotropic diffusion,'' in
  \emph{Image Processing, 2005. ICIP 2005. IEEE International Conference on},
  vol.~1, Sept 2005, pp. I--5--8.

\bibitem{Dey}
N.~Dey, L.~Blanc-Feraud, C.~Zimmer, P.~Roux, Z.~Kam, J.-C. Olivo-Marin, and
  J.~Zerubia, ``Richardson–lucy algorithm with total variation regularization
  for 3d confocal microscope deconvolution,'' \emph{Microscopy Research and
  Technique}, vol.~69, no.~4, pp. 260--266, 2006.

\bibitem{dey:inria-00070726}
N.~Dey, L.~Blanc-F{\'e}raud, C.~Zimmer, P.~Roux, Z.~Kam, J.-C. Olivo-Marin, and
  J.~Zerubia, ``{3D Microscopy Deconvolution using Richardson-Lucy Algorithm
  with Total Variation Regularization},'' Research Report RR-5272, 2004.

\bibitem{Pankajakshan:09}
P.~Pankajakshan, B.~Zhang, L.~Blanc-F\'{e}raud, Z.~Kam, J.-C. Olivo-Marin, and
  J.~Zerubia, ``Blind deconvolution for thin-layered confocal imaging,''
  \emph{Appl. Opt.}, vol.~48, no.~22, pp. 4437--4448, Aug 2009.

\bibitem{Zhang:07}
B.~Zhang, J.~Zerubia, and J.-C. Olivo-Marin, ``Gaussian approximations of
  fluorescence microscope point-spread function models,'' \emph{Appl. Opt.},
  vol.~46, no.~10, pp. 1819--1829, Apr 2007.

\bibitem{Com04}
P.~L. Combettes and J.~Pesquet, ``Image restoration subject to a total
  variation constraint,'' \emph{IEEE Transactions on Image Processing},
  vol.~13, pp. 1213--1222, 2004.

\bibitem{Com11}
P.~L. Combettes and J.-C. Pesquet, ``\BIBforeignlanguage{English}{Proximal
  splitting methods in signal processing},'' in
  \emph{\BIBforeignlanguage{English}{Fixed-Point Algorithms for Inverse
  Problems in Science and Engineering}}, ser. Springer Optimization and Its
  Applications, H.~H. Bauschke, R.~S. Burachik, P.~L. Combettes, V.~Elser,
  D.~R. Luke, and H.~Wolkowicz, Eds.\hskip 1em plus 0.5em minus 0.4em\relax
  Springer New York, 2011, pp. 185--212.

\bibitem{Kos12}
K.~Kose, V.~Cevher, and A.~E. Cetin, ``Filtered variation method for denoising
  and sparse signal processing,'' \emph{IEEE International Conference on
  Acoustics, Speech and Signal Processing (ICASSP)}, pp. 3329--3332, 2012.

\bibitem{PesquetICASSP}
G.~Chierchia, N.~Pustelnik, J.-C. Pesquet, and B.~Pesquet-Popescu, ``An
  epigraphical convex optimization approach for multicomponent image
  restoration using non-local structure tensor,'' in \emph{Acoustics, Speech
  and Signal Processing (ICASSP), 2013 IEEE International Conference on}, 2013,
  pp. 1359--1363.

\bibitem{ChierchiaSIVP14}
------, ``\BIBforeignlanguage{English}{Epigraphical projection and proximal
  tools for solving constrained convex optimization problems},''
  \emph{\BIBforeignlanguage{English}{Signal, Image and Video Processing}}, pp.
  1--13, 2014.

\bibitem{Oppenheim80}
M.~Hayes, J.~Lim, and A.~Oppenheim, ``Signal reconstruction from phase or
  magnitude,'' \emph{Acoustics, Speech and Signal Processing, IEEE Transactions
  on}, vol.~28, no.~6, pp. 672--680, Dec 1980.

\bibitem{Oppenheim81}
A.~Oppenheim and J.~Lim, ``The importance of phase in signals,''
  \emph{Proceedings of the IEEE}, vol.~69, no.~5, pp. 529--541, May 1981.

\bibitem{Oppenheim82}
A.~V. Oppenheim, M.~H. Hayes, and J.~S. Lim, ``Iterative procedures for signal
  reconstruction from fourier transform phase,'' \emph{Optical Engineering},
  vol.~21, no.~1, pp. 211\,122--211\,122--, 1982.

\bibitem{Oppenheim84}
S.~Curtis, J.~Lim, and A.~Oppenheim, ``Signal reconstruction from one bit of
  fourier transform phase,'' in \emph{Acoustics, Speech, and Signal Processing,
  IEEE International Conference on ICASSP '84.}, vol.~9, Mar 1984, pp.
  487--490.

\bibitem{Rud92}
L.~I. Rudin, S.~Osher, and E.~Fatemi, ``Nonlinear total variation based noise
  removal algorithms,'' \emph{Physica D: Nonlinear Phenomena}, vol.~60, no.
  1–4, pp. 259 -- 268, 1992.

\bibitem{PesquetTIP}
N.~Pustelnik, C.~Chaux, and J.~Pesquet, ``Parallel proximal algorithm for image
  restoration using hybrid regularization,'' \emph{IEEE Transactions on Image
  Processing}, vol.~20, no.~9, pp. 2450--2462, 2011.

\bibitem{Chambolle}
A.~Chambolle, ``An algorithm for total variation minimization and
  applications,'' \emph{Journal of Mathematical Imaging and Vision}, vol.~20,
  no. 1-2, pp. 89--97, Jan. 2004.

\bibitem{TofighiICIP}
M.~Tofighi, K.~Kose, and A.~E. Cetin, ``Denoising using projections onto the
  epigraph set of convex cost functions,'' in \emph{Image Processing (ICIP),
  2014 IEEE International Conference on}, Oct 2014, pp. 2709--2713.

\bibitem{Cen12}
Y.~Censor, W.~Chen, P.~L. Combettes, R.~Davidi, and G.~Herman,
  ``\BIBforeignlanguage{English}{On the {E}ffectiveness of {P}rojection
  {M}ethods for {C}onvex {F}easibility {P}roblems with {L}inear {I}nequality
  {C}onstraints},'' \emph{\BIBforeignlanguage{English}{Computational
  Optimization and Applications}}, vol.~51, no.~3, pp. 1065--1088, 2012.

\bibitem{GlobalSIP2013}
A.~E. Cetin, A.~Bozkurt, O.~Gunay, Y.~H. Habiboglu, K.~Kose, I.~Onaran, R.~A.
  Sevimli, and M.~Tofighi, ``Projections onto convex sets (pocs) based
  optimization by lifting,'' \emph{IEEE GlobalSIP, Austin, Texas, USA}, 2013.

\bibitem{You82}
D.~Youla and H.~Webb, ``Image restoration by the method of convex projections:
  Part 1 theory,'' \emph{Medical Imaging, IEEE Transactions on}, vol.~1, no.~2,
  pp. 81--94, 1982.

\bibitem{Cetin1990}
A.~E. Cetin and A.~Tekalp, ``Robust reduced update kalman filtering,''
  \emph{Circuits and Systems, IEEE Transactions on}, vol.~37, no.~1, pp.
  155--156, Jan 1990.

\bibitem{Cetin88}
A.~E. \c{C}etin and R.~Ansari, ``Convolution-based framework for signal
  recovery and applications,'' \emph{J. Opt. Soc. Am. A}, vol.~5, no.~8, pp.
  1193--1200, Aug 1988.

\bibitem{KivancSPM}
K.~Kose and A.~Cetin, ``Low-pass filtering of irregularly sampled signals using
  a set theoretic framework [lecture notes],'' \emph{Signal Processing
  Magazine, IEEE}, vol.~28, no.~4, pp. 117--121, July 2011.

\bibitem{Tru85}
H.~Trussell and M.~R. Civanlar, ``The {L}andweber {I}teration and {P}rojection
  {O}nto {C}onvex {S}et,'' \emph{IEEE Transactions on Acoustics, Speech and
  Signal Processing}, vol.~33, no.~6, pp. 1632--1634, 1985.

\bibitem{Stark:81}
H.~Stark, D.~Cahana, and H.~Webb, ``Restoration of arbitrary finite-energy
  optical objects from limited spatial and spectral information,'' \emph{J.
  Opt. Soc. Am.}, vol.~71, no.~6, pp. 635--642, Jun 1981.

\bibitem{Combettes93}
P.~Combettes, ``The foundations of set theoretic estimation,''
  \emph{Proceedings of the IEEE}, vol.~81, no.~2, pp. 182--208, Feb 1993.

\bibitem{Ayers:88}
G.~R. Ayers and J.~C. Dainty, ``Iterative blind deconvolution method and its
  applications,'' \emph{Opt. Lett.}, vol.~13, no.~7, pp. 547--549, Jul 1988.

\bibitem{Slavakis}
S.~Theodoridis, K.~Slavakis, and I.~Yamada, ``Adaptive learning in a world of
  projections,'' \emph{Signal Processing Magazine, IEEE}, vol.~28, no.~1, pp.
  97--123, Jan 2011.

\bibitem{Yamada2001473}
I.~Yamada, ``The hybrid steepest descent method for the variational inequality
  problem over the intersection of fixed point sets of nonexpansive mappings,''
  in \emph{Inherently Parallel Algorithms in Feasibility and Optimization and
  their Applications}, ser. Studies in Computational Mathematics, Y.~C.
  Dan~Butnariu and S.~Reich, Eds.\hskip 1em plus 0.5em minus 0.4em\relax
  Elsevier, 2001, vol.~8, pp. 473 -- 504.

\bibitem{Vonesch09}
C.~Vonesch and M.~Unser, ``A fast multilevel algorithm for wavelet-regularized
  image restoration,'' \emph{Image Processing, IEEE Transactions on}, vol.~18,
  no.~3, pp. 509--523, March 2009.

\end{thebibliography}

\end{document}